\documentclass{amsart}

\usepackage{amsmath,amsthm,amsfonts,amssymb,mathrsfs,mathtools}
\usepackage[
top=1in, left=1.5in, right=1in, bottom=1in]{geometry}
\usepackage{setspace}
\usepackage[nice]{nicefrac}
\usepackage[all,cmtip]{xy}
\usepackage{multind}
\usepackage{lscape}
\usepackage{varioref}
\usepackage{cite}
\usepackage{fancyhdr}

\makeindex{Index}
\makeindex{Notation}

\xyoption{arc}

\DeclareMathAlphabet{\mathpzc}{OT1}{pzc}{m}{it}

\newcommand{\sgn}{\operatorname{sgn}}

\newcommand{\C}{\mathbb{C}}

\newcommand{\Z}{\mathbb{Z}}

\newcommand{\K}{\mathpzc{K}}

\newcommand{\pK}[2]{\mathpzc{K}({#1},{#2})}
\newcommand{\ksig}[1]{K_{_{\Sigma_{#1}}}(X^{#1})}

\newcommand{\kwr}[2]{K_{_{\Sigma_{#1}\wr\Sigma_{#2}}}(X^{#1#2})}

\newcommand{\Ind}[2]{\operatorname{Ind}_{_{#1}}^{^{#2}}}
\newcommand{\Res}[2]{\operatorname{Res}_{_{#1}}^{^{#2}}}
\newcommand{\Vect}{\operatorname{Vect}}
\newcommand{\pt}{\operatorname{pt}}

\theoremstyle{definition}
\newtheorem{theorem}{Theorem}

\newtheorem{definition}{Definition}[section]

\newtheorem{lemma}[theorem]{Lemma}
\newtheorem*{definition*}{Definition}

\title{$\K$: An Equivariant $K$-Theory Functor from Spaces to $\lambda$-rings}

\author{Joseph C. Johnson}

\email{jjohnson@mbc.edu}

\address{Mathematics Department,
         Mary Baldwin College,
         Staunton, Virginia, 24401,
         United States}

\begin{document}

\begin{abstract}

For a compact Hausdorf space $X$, let
$$
\K(X)=\bigoplus_{n\geq 0}K_{_{\Sigma_n}}(X^n).
$$
We show that $\K$ takes values in $\lambda$-rings and satisfies a Thom isomorphism.  In the case that $X$ is constructed entirely of even-dimensional cells, $\K(X)$ is the free $\lambda$-ring on generators in one-to-one correspondence with those cells.

\end{abstract}

\maketitle

\tableofcontents

\section{Introduction}
Grothendieck first defined $\lambda$-rings in some unpublished work in 1957 and published in \cite{Grothendieck} under the name ``special $\lambda$-rings''.  In short a $\lambda$-ring is a ring $A$, with certain set functions $\lambda^n:A\longrightarrow A$ that satisfy the relations that are satisfied by the exterior powers of vector spaces.  Examples of $\lambda$-rings include complex representation rings of groups and K-theory groups of spaces.  In both cases the $\lambda$-ring functions come from exterior power of vector spaces.  

In this paper, we introduce a homotopy functor from suitable pairs of topological spaces to $\lambda$-rings.  For a compact, Hausdorf $G$-space $X$ and a closed subspace $A$, let $K_G(X,A)$ be the $G$-equivariant K-theory of the pair $(X,A)$.  The notation $(X,A)^n$ denotes the pair
$$
(X^n,(X^{n-1}\times A)\cup (X^{n-2}\times A\times X)\cup\ldots\cup (A\times X^{n-1})).
$$
These are $\Sigma_n$-spaces.  We then let
$$
\K(X,A)=\bigoplus_{n\geq 0}K_{_{\Sigma_n}}((X,A)^n)
$$
and write $\K(X,\emptyset)$ as $\K(X)$.

\subsection{Main Results}\text{}
\\
Our first result is
\begin{theorem}\label{KisRing}
$\K(X,A)$ is a graded $\lambda$-ring.  Furthermore, if each $K_{_{\Sigma_n}}((X,A)^n)$ is a finitely generated and free and $K^1_{_{\Sigma_n}}((X,A)^n)=0$ for all $n$, then $\pK{X}{A}$ is also equipped with a coproduct that makes it into a Hopf algebra.
\end{theorem}

The structure on $\K(X,A)$ arises as follows.  For each $i,j\geq 0$ we have an inclusion $\Sigma_i\times\Sigma_j\leq \Sigma_{i+j}$.  This induces maps
$$
K_{_{\Sigma_i}}((X,A)^i) \otimes K_{_{\Sigma_j}}((X,A)^j)\stackrel{\times}{\longrightarrow} K_{_{\Sigma_i\times\Sigma_j}}((X,A)^{i+j})\stackrel{\Ind{\Sigma_i\times\Sigma_j}{\Sigma_{i+j}}}{\longrightarrow} K_{_{\Sigma_{i+j}}}((X,A)^{i+j})
$$
and
$$
K_{_{\Sigma_{i+j}}}((X,A)^{i+j})\stackrel{\Res{\Sigma_i\times\Sigma_j}{\Sigma_{i+j}}}{\longrightarrow} K_{_{\Sigma_i\times\Sigma_j}}((X,A)^{i+j}) \stackrel{\times}{\longleftarrow} K_{_{\Sigma_i}}((X,A)^i)\otimes K_{_{\Sigma_j}}((X,A)^j).
$$
These give the product, and in the cases listed above, the coproduct on $\K(X,A)$.

It is also good to note that the functor $\K$ satisfies excision.  A functor of pairs, $F$, satisfies excision if when $X=A\cup B$ is the union of two closed subspaces, then
$$
F(X,A)\cong F(B,B\cap A).
$$
So, suppose that $X=A\cup B$.  Then
$$
X^n=B^n\cup[(X^{n-1}\times A) \cup\ldots\cup (A\times X^{n-1})]
$$ 
and
$$
B^n\cap [(X^{n-1}\times A)\cup\ldots\cup (A\times X^{n-1})]=[B^{n-1}\times (A\cap B)]\cup \ldots\cup [(A\cap B)\times B^{n-1})].
$$
Since each $K_{_{\Sigma_n}}$ is a cohomology functor, it satisfies excision.  Using the above, we have
$$
K_{_{\Sigma_n}}((X,A)^n)\cong K_{_{\Sigma_n}}((B,B\cap A)^n).
$$
Adding these up we get

\begin{theorem}(Excision)\label{Excision}
If $X=A\cup B$ is a union of two closed subspaces, then
$$
\K(X,A)\cong \K(B,B\cap A).
$$
\end{theorem}

We obtain a $\lambda$-ring structure as follows.  If $E\in K_{_{\Sigma_i}}((X,A)^i)$, the $n$-fold exterior product $E^{\boxtimes n}\in K_{_{\Sigma_i\wr\Sigma_n}}((X,A)^{in})$.  Then $\lambda^n(E)\in K_{_{in}}((X,A)^{in})$ is defined by
$$
\lambda^{n}(E)=\Ind{\Sigma_i\wr\Sigma_n}{\Sigma_{in}}(\sgn_n\otimes E^{\boxtimes n})
$$
where $\sgn_n$ is the sign representation of $\Sigma_n$.

Our next main result is that the functor $\K$ satisfies a Thom isomorphism.

\begin{theorem}\label{Thom}
If $E\longrightarrow X$ is a vector bundle with associated disk and sphere bundles $D(E)$ and $S(E)$, there is an isomorphism of $\lambda$-rings:
$$
\K(X,A)\longrightarrow \K(D(E),S(E)\cup D(E)|_A).
$$
\end{theorem}

Our final result concerns the restriction of $\K$ to a certain subclass of spaces.  Let $\mathfrak{F}$ be the free $\lambda$-ring functor from the category of abelian groups to $\lambda$-rings.

\begin{theorem}\label{Free}
If $X$ is obtained from $A$ by attaching finitely many even-dimensional cells, then the natural map of $\lambda$-rings
$$
\mathfrak{F}(K(X,A))\longrightarrow \K(X,A)
$$
is an isomorphism.
\end{theorem}

\subsection{Relation to the Literature}\text{}
\\
This functor $\K$ is not new.  We remark that 
$$
\K(\operatorname{pt})=\bigoplus_{n\geq 0}R(\Sigma_n)
$$
has been known to be the free $\lambda$-ring on one generator at least since \cite{MacDonald} and our work presents another proof of this fact.

In response to a comment by Grojnowski \cite{Grojnowski}, Graeme Segal, studied the object
$$
\bigoplus_{n\geq 0}K_{_{\Sigma_n}}(X^n)\otimes\C
$$
in an unpublished note \cite{UnSegal}.  Weiqiang Wang replaced $\Sigma_n$ with $G\wr\Sigma_n$, the wreath product of $\Sigma_n$ with a finite group $G$ \cite{Wang}.  He proves that
$$
\bigoplus_{n\geq 0}K_{_{G\wr\Sigma_n}}(X^n)\otimes\C
$$
is the free lambda ring on $K_G(X)$.  Moreover, it is shown that $\K(X)\otimes\C$ is a Fock space associated to a certain Heisenberg algebra.  

Another connection to the literature is as follows:  There are maps
$$
K_{_{\Sigma_n}}((X,A)^n) \longrightarrow K(E\Sigma_n \times_{\Sigma_n} (X,A)^n)
$$
which are algebraic completions by the Atiyah-Segal Completion Theorem \cite[Proposition 4.2]{SegalCompletion}.  Let
$$
\mathbb{P}(X,A)=\coprod_{n\geq 0} E\Sigma_n\times_{\Sigma_n}(X,A)^n.
$$
Combined, the above maps define a natural map
$$
\K(X,A)\longrightarrow K(\mathbb{P}(X,A)).
$$
Furthermore, if $X$ is path connected, $K(\mathbb{P}(X,\operatorname{pt}))$ identifies with $K(\Omega^\infty \Sigma^\infty X)$.  This isomorphism will not necessarily be an isomorphism of $\lambda$-rings.  In general the former will be an ``associated graded" version of the latter, so that $K(\Omega^\infty \Sigma^\infty X)$ can be calculated as an appropriate algebraic completion of $\K(X,\operatorname{pt})$.

In 2011, $\K(X)$ has shown up in the work of Nora Ganter and Mikhail Kapranov as the Grothendieck group  of the symmetric powers in the category of coherent sheaves on a smooth projective variety \cite{Ganter}.

\subsection{Methods of Proof}\text{}
\\
We prove that $\K(X,A)$ is a $\lambda$-ring by showing that it is a $\tau$-ring in the sense of \cite{Hoffman} and use his result that every $\tau$-ring has the structure of a $\lambda$-ring.  The categories are, in fact, equivalent.  We do this in Section \ref{KisTau}, after a brief introduction to $\lambda$-rings and $\tau$-rings in Section \ref{LambdaandTau}.  We prove Theorem \ref{Thom} in Section \ref{ThomIso} and Theorem \ref{Free} in Section \ref{FreeIso}.  We prove Theorem \ref{Thom} by showing that the usual Thom isomorphism can be extended to the structure on $\K(X,A)$.  We prove Theorem \ref{Free} by induction on cells.

\section{$\lambda$-Rings and $\tau$-Rings}\label{LambdaandTau}

We begin with the definition of a $\lambda$-ring.

\begin{definition}
A $\lambda$-ring is a commutative ring with unit, $A$, and a set of functions for each $n=0,1,2,\ldots$
$$
\lambda^n: A\longrightarrow A
$$
satisfying:
\begin{itemize}
\item[(i)] $\lambda^0(x)=1$, $\lambda^1(x)=x$ for $x \in A$
\item[(ii)] $\lambda^k(x+y)=\sum_{i=0}^k\lambda^i(x)\lambda^{k-i}(y)$ for $x,y\in A$
\item[(iii)]
$\lambda^k(xy)=p_k(\lambda^1(x),\ldots ,\lambda^k(x),\lambda^1(y),\ldots ,\lambda^k(y))$ for $x,y\in A$.  Here $p_k$ is the unique polynomial such that
$$
\prod_{i,j}(1+a_ib_jt)  =\sum_kp_k(\sigma_1(\bar{a}),\ldots ,\sigma_k(\bar{a}),\sigma_1(\bar{b}),\ldots ,\sigma_k(\bar{b}))t^k
$$
The $\sigma_i$ are the elementary symmetric polynomials.
\item[(iv)]

$\lambda^k(\lambda^l(x))=q_{k,l}(\lambda^1(x),\ldots,\lambda^{kl}(x))$ for $x\in A$.  Here $q_{k,l}$ is the unique polynomial such that 
$$
\prod_{i_1<\cdots<i_l}(1+x_{i_1}\cdots x_{i_l}t)=\sum_kq_{k,l}(\sigma_1(\bar{x}),\ldots,\sigma_{kl}(\bar{x}))t^k
$$
in $\Z[x_1,\ldots,x_n,t]$.
\end{itemize}
\end{definition}

A basic example of a $\lambda$-ring is the $G$-equivariant $K$-theory of a compact Hausdorf space, $X$, where the $\lambda$-operations are the exterior powers of vector bundles.  

Now for the definition of a $\tau$-ring.  Suppose that $A$ is a commutative ring with unit.  Define
$$
AP=\prod_{n\geq 0}A\otimes R(\Sigma_n)
$$
where $R(\Sigma_n)$ is the representation ring of the $n$th symmetric group.  Furthermore let
$$
AP^{(2)}=\prod_{n,m\geq 0}A\otimes R(\Sigma_n)\otimes R(\Sigma_m),\quad AQ=\prod_{n,m\geq 0}A\otimes R(\Sigma_m\wr\Sigma_n)
$$  

We have a variety of maps that we must define before we can define a $\tau$-ring.  First of all, we have a multiplication map
$$
\times: AP\otimes AP\longrightarrow AP.
$$
The $n$th coordinate of this map is given by the composite
$$
AP\longrightarrow \prod_{i+j=n} A\otimes R(\Sigma_i)\otimes A\otimes R(\Sigma_j)\longrightarrow A\otimes R(\Sigma_n).
$$
Where the first map is projection and the second map is the induction multiplication in $R_*$. 

Our next multiplication comes from the levelwise internal tensor product on each $R(\Sigma_n)$.  This will induce a map
$$
\cdot:AP\otimes AP\longrightarrow AP
$$
with $n$th coordinate
$$
AP\otimes AP\longrightarrow A\otimes R(\Sigma_n)\otimes A\otimes R(\Sigma_n)\longrightarrow A\otimes R(\Sigma_n)
$$
where the first map is projection and the second sends $a_1\otimes V \otimes  a_2\otimes W$ to $a_1a_2\otimes (V\otimes W)$.

We have a map
$$
\Box: AP\longrightarrow AQ.
$$
This is induced by the restriction maps $R(\Sigma_{mn})\rightarrow R(\Sigma_m\wr \Sigma_n)$ and has $(n,m)$ coordinate function:
$$
AP\longrightarrow A\otimes R(\Sigma_{mn})\longrightarrow A\otimes R(\Sigma_m\wr\Sigma_n).
$$
Where the first map is projection and the second is the aforementioned restriction.  

From the coproduct structure on $\oplus R(\Sigma_n)$ we obtain a map 
$$
\Delta: AP\rightarrow AP^{(2)},
$$
with $(i,j)$ coordinate function
$$
AP\longrightarrow A\otimes R(\Sigma_n)\longrightarrow A\otimes R(\Sigma_i)\otimes R(\Sigma_j).
$$
Where $n=i+j$; the first map is projection; and the second map is the coproduct.  

Next, we have a map
$$
\mu: AP\otimes AP\longrightarrow AP^{(2)}.
$$
The $(i,j)$ coordinate is the composite
$$
AP\otimes AP \longrightarrow A\otimes R(\Sigma_i)\otimes A\otimes R(\Sigma_j)\longrightarrow A\otimes R(\Sigma_i)\otimes R(\Sigma_j).
$$
Where the first map is projection and the second is multiplication in $A$.  

Finally, let $e_n$ denote the trivial representation of $R(\Sigma_n)$, with the convention that $e_0=1\in\Z$.  With this notation we may now introduce the notion of a $\tau$-ring.
\begin{definition}
A $\tau$-ring is a commutative ring with unit, $A$, and a function 
$$
\tau: A\rightarrow AP
$$ 
satisfying:
\begin{enumerate}
\item  $\tau(x)\in 1\otimes e_0 + x\otimes e_1 +\prod_{n\geq 2}A\otimes R(\Sigma_n)$ for all $x\in A$.
\item  $\tau(x+y)=\tau(x)\times\tau(y).$ for $x,y\in A$
\item  The following commutes:
$$
\xymatrix{ A\ar[r]^-{\tau^{(2)}}\ar[d]_\tau & AP\otimes AP\ar[d]^\mu \\
AP\ar[r]^\Delta & AP^{(2)}}
$$
where $\tau^{(2)}$ is the map $a\mapsto \tau(a)\otimes\tau(a)$.

\item  $\tau(xy)=\tau(x)\cdot\tau(y)$ for $x,y\in A$.
\item  The following commutes:
$$
\xymatrix{ A\ar[r]^-{\tau}\ar[d]_{\tau} &AP\ar[d]^{\dot{\tau}} \\
AP\ar[r]^-{\Box} & AQ}
$$
The definition of $\dot{\tau}$ is a little involved:  Assume $\tau: A\rightarrow AP$ is a map that satisfies the first three properties for a $\tau$-ring.  Let $\tau^k$ be the $k$th component of $\tau$.  We have the map
$$
A\times R(\Sigma_l) \longrightarrow A\otimes R(\Sigma_k)\otimes R(\Sigma_l\wr\Sigma_k) 
$$
which takes $(x,V)$ to $\tau^k(x)\otimes V^{\otimes k}$.  Next we have a projection map 
$$
p:\Sigma_l\wr \Sigma_k\rightarrow \Sigma_k
$$ 
that induces a pullback map 
$$
p^*: R(\Sigma_k)\rightarrow R(\Sigma_l \wr \Sigma_k).
$$  
This leads to another map
$$
\operatorname{Id}_A\otimes p^* \otimes \operatorname{Id}_{R(\Sigma_l\wr\Sigma_k)}:A\otimes R(\Sigma_k)\otimes R(\Sigma_l\wr\Sigma_k)\rightarrow A\otimes R(\Sigma_l\wr\Sigma_k)\otimes R(\Sigma_l\wr \Sigma_k).
$$
Finally, we note that there is an internal tensor product on $R(\Sigma_l\wr\Sigma_k)$, call it $\gamma$.  Our last map is:
$$
\operatorname{Id}_A\otimes \gamma:A\otimes R(\Sigma_l\wr\Sigma_k)\otimes R(\Sigma_l\wr \Sigma_k) \longrightarrow A\otimes R(\Sigma_l\wr\Sigma_k).
$$
The composition of all the above maps is denoted 
$$
\dot{\tau}^k_l:A\times R(\Sigma_l)\longrightarrow A\otimes R(\Sigma_l\wr \Sigma_k).
$$
We define $\dot{\tau}_l$ to be the product over all $k$.  Since $\tau$ satisfies the first three properties of a $\tau$-ring, it is bi-additive from addition to cross product.  Thus it extends uniquely to a function
$$
\dot{\tau}_l: A\otimes R(\Sigma_l)\longrightarrow A\otimes \prod_{k\geq 0}R(\Sigma_l\wr\Sigma_k).
$$
Put 
$$\dot{\tau}=\prod \dot{\tau}_l.
$$ 
\end{enumerate}
\end{definition}

With our definition of a $\tau$-ring in place, let us now show how to obtain a $\lambda$-ring from a $\tau$-ring.  For $n=1,2,3,\dots$ let
$$
g_n: AP\longrightarrow A
$$
be the composite
$$
AP\longrightarrow A\otimes R(\Sigma_n)\stackrel{g_n'}{\longrightarrow} A.
$$
Here the first map is projection and the second is
\begin{equation*}
g_n'(a\otimes \nu_i) =
\begin{cases}
 a & \text{if $\nu_i=\sgn_n$}\\
 0 & \text{otherwise}
\end{cases}
\end{equation*}
Here the $\nu_i$ are the irreducible representations of $\Sigma_n$ and $\sgn_n$ is the sign representation of $\Sigma_n$.  We may define $\lambda$ operations by $\lambda^n(a)=g_n(\tau(a))$.

\section{$\K(X,A)$ is a $\tau$-ring}\label{KisTau}

The proof presented is only for $\K(X)$, but is easily seen to be valid for the functor of pairs $\pK{X}{A}$.  It is presented in this way to provide ease of exposition.  

Let us introduce some notation.  Let
$$
\K P=\prod_{l\geq 0}\K(X)\otimes R(\Sigma_l),
$$
$$
\K P^{(2)}=\prod_{i,j\geq 0}\K(X)\otimes R(\Sigma_i)\otimes R(\Sigma_j),
$$
and
$$
\K Q=\prod_{k,l\geq 0}\K(X)\otimes R(\Sigma_l\wr\Sigma_k).
$$
Maps $\tau^m_n$ will be defined for $n,m\geq 1$.  Then, $\tau$ will simply be the product over all $n$ and $m$.

Let $E$ be a $\Sigma_n$-bundle over $X^n$.  First we have the $m$th power map for $m\geq 1$:
$$
\mathcal{P}^m_n: \ksig{n}\rightarrow\kwr{n}{m}
$$
that sends $E$ to $E^{\boxtimes m}$, the external tensor product.  Here $\Sigma_m$ acts by $g\cdot (v_i)=(v_{g^{-1}(i)})$ and $\Sigma_n$ acts diagonally.  The power map and some of its elementary properties are discussed in \cite[Section 1]{Rezk}.  A second family of maps is constructed as follows.  Embed $\Sigma_n\wr\Sigma_m$ as a subgroup of $\Sigma_m\times\Sigma_{nm}$, via the product of the canonical quotient map $\Sigma_n\wr\Sigma_m\rightarrow\Sigma_m$ with the usual inclusion map $\Sigma_n\wr\Sigma_m<\Sigma_{mn}$ given by concatenation.  If we regard $X^{mn}$ as a $\Sigma_m\times\Sigma_{mn}$ space where the first factor acts trivially and the second acts in the standard manner, the restriction of this action to our embedded subgroup is the standard action of $\Sigma_n\wr\Sigma_m$ on $X^{mn}$.  We can therefore define the map for $m\geq 1$:
$$
\delta^m_n= \operatorname{Ind}_{\Sigma_n\wr\Sigma_m}^{\Sigma_m\times\Sigma_{mn}}:\kwr{n}{m}\longrightarrow K_{_{\Sigma_m\times\Sigma_{mn}}}(X^{mn})\cong\ksig{mn}\otimes R(\Sigma_m).
$$
The induction map is also discussed in \cite[Section 1]{Rezk} and is a generalization of the induction in the groups setting.  We then define $\tau_n^m$ to be the composition of $\mathcal{P}^m_n$ followed by $\delta^m_n$.  The $\tau$ map is then the product over all $m$ and $n$ of $\tau_n^m$.

Before we begin the proof, it should be noted that the resulting $\lambda$-ring structure is simply
$$
\lambda^m(E)=\Ind{\Sigma_n\wr\Sigma_m}{\Sigma_{nm}}(\sgn_m\otimes E^{\boxtimes m}).
$$
In the case that $X=\operatorname{pt}$, this is the plethysm operation.

In order to show that this is indeed a $\tau$-ring structure a few things must be shown.  Each will be introduced and proven in turn.  First we show that it is a pre-$\tau$-ring, i.e. it satisfies these first three properties:
\begin{itemize}

\item[1.] $\tau(x)\in 1\otimes e_0+x\otimes e_1 +\prod_{l\geq 2}\K(X)\otimes R(\Sigma_l)$, where $e_0$ and $e_1$ are the trivial representations of $\Sigma_0$ and $\Sigma_1$.
\begin{proof}
First, we have assumed that $\tau^0_n$ is the constant map $1\otimes e_0$.  Now, the map $\mathcal{P}^1_n$ is just the identity map on $\ksig{n}$.  And for $m=1$ the embedding of $\Sigma_n=\Sigma_n\wr\Sigma_1$ in $\Sigma_n=\Sigma_1\times\Sigma_n$ is again the identity and so $x\in\ksig{n}$ certainly maps to $x\otimes e_1$.
\end{proof}

\item[2.] $\tau(x+y)=\tau(x)\times\tau(y)$
\begin{proof}
Notice that we already have a cross product structure on $\prod_{l\geq 0}\kwr{n}{l}$ given by the composition of:
$$
\kwr{n}{m}\otimes\kwr{n}{k} \stackrel{\boxtimes}{\longrightarrow}K_{_{\Sigma_n\wr\Sigma_m\times\Sigma_n\wr\Sigma_k}}(X^{nm}\times X^{nk})
$$
with
$$
K_{_{\Sigma_n\wr\Sigma_m\times\Sigma_n\wr\Sigma_k}}(X^{nm}\times X^{nk})\stackrel{\operatorname{Ind}}{\longrightarrow}K_{_{\Sigma_n\wr\Sigma_{m+k}}}(X^{n(k+m)}).
$$
Since $\delta^m_n$ is a ring homomorphism with the appropriate multiplications we need only show the claim for $\mathcal{P}^m_n$.  For $m=0$, the equation
$$
\mathcal{P}^0(x+y)=\mathcal{P}^0(x)\times\mathcal{P}^0(y)
$$
is trivial.  For $m>0$, it follows from the isomorphism
$$
(E_0\oplus E_1)^{\boxtimes m}\cong\oplus_{i=0}^m W_i
$$
of $\Sigma_n\wr\Sigma_m$-bundles over $X^{nm}$, where $W_i=\bigoplus_{J_i}E_{i_1}\boxtimes\cdots\boxtimes E_{i_m}$, and the index set, $J_i$, is
$$
\{(i_1,\ldots,i_m)\mid i_l=\text{0 or 1}, \sum_l i_l=i\}.
$$
Each $W_i$ is a direct sum of $\begin{pmatrix} m\\i \end{pmatrix}$ sub-bundles which are permuted transitively by $\Sigma_n\wr\Sigma_m$.  Also one of the summands is $E_0^{\boxtimes m-i}\boxtimes E_1^{\boxtimes i}$ which is invariant under the inclusion of $\Sigma_n\wr\Sigma_{m-i}\times\Sigma_n\wr\Sigma_i$ in $\Sigma_n\wr\Sigma_m$.  Therefore $W_i\cong\mathcal{P}^{m-i}(E_0)\times\mathcal{P}^i(E_1)$.  The result then follows.
\end{proof}

\item[3.] The following commutes:
$$
\xymatrix{ \K(X)\ar[r]^-{\tau^{(2)}}\ar[d]_\tau & \K P\otimes \K P\ar[d]^\mu \\
\K P\ar[r]^\Delta & \K P^{(2)}}
$$
where $\mu$ is multiplication in the left hand factor, $\K(X)$, of $\K P$ and $\Delta$ is induced by the restrictions maps from $\Sigma_m$ to $\Sigma_{m-i}\times\Sigma_i$.
\begin{proof}
To prove this we shall split the diagram into two easier pieces:
$$
\xymatrix{ \ksig{n} \ar[r]^-{\mathcal{P}^{(2)}} \ar[d]_{\mathcal{P}} &\prod_{i,j}K_{_{\Sigma_n\wr\Sigma_i\times\Sigma_n\wr\Sigma_j}}(X^{ni}\times X^{nj})\ar[r]^-{\delta^{(2)}} & \K P\otimes \K P\ar[d]^{\mu}\\
\prod_l K_{_{\Sigma_n\wr\Sigma_l}}(X^{nl})\ar[ru]^{\Delta} \ar[r]^{\delta} & \K P\ar[r]^{\Delta} & \K P^{(2)}}
$$
where the $\Delta$ maps are the ``coproducts'' induced by the restrictions to the subgroups $\Sigma_i\times\Sigma_j\leq \Sigma_{i+j}$ and $\Sigma_n\wr\Sigma_i\times\Sigma_n\wr\Sigma_j\leq\Sigma_n\wr\Sigma_{i+j}$.  The left triangle commutes due to the fact that the canonical isomorphism $E^{\boxtimes (i+j)}\cong E^{\boxtimes i}\boxtimes E^{\boxtimes j}$ is an isomorphism of $\Sigma_n\wr\Sigma_i\times\Sigma_n\wr\Sigma_j$-bundles.

Now suppose $E$ is a $\Sigma_n\wr\Sigma_l$-bundle over $X^{nl}$.  Then
$$
\Delta(E)=\bigoplus_{i+j=l}\Res{\Sigma_n\wr\Sigma_i\times\Sigma_n\wr\Sigma_j}{\Sigma_n\wr\Sigma_l}(E).
$$
Then $\delta^{(2)}$ followed by $\mu$ is just
$$
\bigoplus_{i+j=l}\Ind{\Sigma_n\wr\Sigma_i\times\Sigma_n\wr\Sigma_j}{\Sigma_{nl}\times\Sigma_i\times\Sigma_j}(\Res{\Sigma_n\wr\Sigma_i\times\Sigma_n\wr\Sigma_j}{\Sigma_n\wr\Sigma_l}(E)).
$$
On the other hand
$$
\delta(E)=\Ind{\Sigma_n\wr\Sigma_l}{\Sigma_{nl}\times\Sigma_l}(E)
$$
and 
$$
\Delta(\delta(E))=\bigoplus_{i+j=l}\Res{\Sigma_{nl}\times\Sigma_i\times\Sigma_j}{\Sigma_l\times\Sigma_{nl}}\Ind{\Sigma_n\wr\Sigma_l}{\Sigma_{nl}\times\Sigma_l}(E).
$$
Conveniently
$$
(\Sigma_n\wr\Sigma_l)\cap(\Sigma_{nl}\times\Sigma_i\times\Sigma_j)=\Sigma_n\wr\Sigma_i\times\Sigma_n\wr\Sigma_j
$$
in $\Sigma_l\times\Sigma_{nl}$.  So, an application of Mackey's Lemma will tell us that these two bundles are isomorphic.
 
\end{proof}

There are two more things to prove in order to show that $\K(X)$ is a $\tau$-ring.  First let's record certain lemmas about wreath products that will be useful.  Note that if there is a group homomorphism $H\rightarrow G$ we get an induced map $H\wr\Sigma_n\rightarrow G\wr\Sigma_n$.  

\begin{lemma}\label{C1}
Let $H$ be a subgroup of a finite group $G$, and $X$ be a $G$-space.  The following commutes:
$$
\xymatrix{ K_H(X)\ar[r]^-{\mathcal{P}^k}\ar[d]_{\operatorname{Ind}} & K_{H\wr\Sigma_k}(X^k)\ar[d]^{\operatorname{Ind}}\\
K_G(X)\ar[r]^-{\mathcal{P}^k} & K_{G\wr\Sigma_k}(X^k)}
$$
\begin{proof}
If $E\in\Vect_H(X)\cong\Vect_G(G\times_H X)$, then $\Ind{H}{G}(E)$ is a product of $E$'s.  Let $s_1,\ldots, s_r$ be a list of coset representatives.  Then
$$
\Ind{H}{G}(E)=\bigoplus_{i=1}^r E_{s_i}.
$$
An element $g\in G$ acts on a vector $(e_{s_1},\ldots,e_{s_r})$ diagonally 
$$
g\cdot (e_{s_1},\ldots,e_{s_r})=(ge_{g^{-1}s_1},\ldots,ge_{g^{-1}s_r}).
$$
We then apply $\mathcal{P}^k$.  Then $(g_1,\ldots, g_k,\sigma)\in G\wr\Sigma_k$ acts by sending a vector
$$
(e_{s_{1_1}},\ldots,e_{s_{r_1}})\boxtimes\cdots\boxtimes (e_{s_{1_k}},\ldots,e_{s_{r_k}})  
$$
to
$$
(g_1e_{g_1^{-1}s_{1_{\sigma^{-1}(1)}}},\ldots, g_ke_{g_k^{-1}s_{r_{\sigma^{-1}(1)}}})\boxtimes\cdots \boxtimes(g_ke_{g_k^{-1}s_{1_{\sigma^{-1}(k)}}},\ldots, g_ke_{g_k^{-1}s_{r_{\sigma^{-1}(k)}}}).
$$
In order for a vector 
$$
\boxtimes_{i=1}^k(e_{s_{1_i}},\ldots,e_{s_{r_i}})
$$
to be invariant under the action of $H\wr\Sigma_k$, we must have $e_{s_{j_i}}=e_{s_{k_i}}$ for all $j$ and $k$.  But, this is isomorphic to $\mathcal{P}^k(E)$.
\end{proof}
\end{lemma}

\begin{lemma}\label{L2}
 The following diagram commutes:
$$
\xymatrix{ \ksig{n}\ar[r]^-{\mathcal{P}^l}\ar[d]_{\mathcal{P}^{kl}} & \kwr{n}{l} \ar[d]^{\mathcal{P}^k}\\
\kwr{n}{kl}\ar[r]^-{\operatorname{Res}} & K_{_{(\Sigma_n\wr\Sigma_l)\wr\Sigma_k}}(X^{kl}) }
$$
\begin{proof}
 The isomorphism $(E^{\boxtimes l})^{\boxtimes k}\cong E^{\boxtimes kl}$ is a map of $(\Sigma_n\wr\Sigma_l)\wr\Sigma_k$-bundles.
\end{proof}

\end{lemma}

Now that we have these lemmas in place, we are equipped to show that $\K(X)$ satisfies the final two properties.
\item[4.]$\tau(xy)=\tau(x)\cdot\tau(y)$
\begin{proof}
 We must show that the following diagram commutes:
$$
\xymatrix{ \ksig{i}\otimes\ksig{j} \ar[d]^-{\mathcal{P}^k\otimes\mathcal{P}^k} \ar[r]^\times & \ksig{i+j} \ar[d]^{\mathcal{P}^k} \\ 
\kwr{i}{k}\otimes\kwr{j}{k} \ar[r]^{*} \ar[d]^-{\delta^k\otimes\delta^k} & \kwr{i+j}{k}\ar[d]^{\delta^k} \\
\ksig{ki}\otimes R(\Sigma_k) \otimes\ksig{kj}\otimes R(\Sigma_k)  \ar[r]^-{\times\otimes\cdot} \ar[d]^-{\Ind{\Sigma_{ki}\times\Sigma_{kj}}{\Sigma_{k(i+j)}}} &  \ksig{k(i+j)}\otimes R(\Sigma_k) \\
\ksig{k(i+j)}\otimes R(\Sigma_k)\otimes R(\Sigma_k) \ar[ru]_{\Res{\Sigma_k}{\Sigma_k\times\Sigma_k}}  }
$$
This new multiplication, $*$, is given by the composition of the isomorphism
$$
\kwr{i}{k}\times\kwr{j}{k}\stackrel{\cong}{\longrightarrow} K_{_{\Sigma_i\wr\Sigma_k\times\Sigma_j\wr\Sigma_k}}(X^{k(i+j)})
$$
with
$$
K_{_{\Sigma_i\wr\Sigma_k\times\Sigma_j\wr\Sigma_k}}(X^{k(i+j)})\stackrel{\operatorname{Res}}{\longrightarrow} K_{_{(\Sigma_i\times\Sigma_j)\wr\Sigma_k}}(X^{k(i+j)})\stackrel{\operatorname{Ind}}{\longrightarrow} \kwr{i+j}{k}.
$$
For the bottom triangle and middle square we again apply Mackey's Lemma to 
$$
(\Sigma_{k(i+j)}\times\Sigma_k)\cap((\Sigma_i\wr\Sigma_k)\times(\Sigma_j\wr\Sigma_k))=(\Sigma_i\times\Sigma_j)\wr\Sigma_k.
$$
\hspace*{6mm} The top square is a little more involved.  The map $\mathcal{P}^k(x\times y)$ is bi-additive in the variables $x$ and $y$ into the cross product invertibles, $1+\prod_{k\geq 1}K_{\Sigma_{i+j}\wr\Sigma_k}$.  Also,
\begin{align*}
\mathcal{P}^k(x+y)*\mathcal{P}^k(z) &= \left(\sum_{i=0}^k\mathcal{P}^i(x)\times\mathcal{P}^{k-i}(y)\right) * \mathcal{P}^k(z)\\
&= \sum_{i=0}^k(\mathcal{P}^i(x)*\mathcal{P}^i(z))\times (\mathcal{P}^{k-i}(y)*\mathcal{P}^{k-i}(z)) 
\end{align*}
So $\mathcal{P}^k(x)*\mathcal{P}^k(y)$ is bi-additive.  This bi-additivity property allows us to divide up the square into:
$$
\xymatrix{ \Vect_{\Sigma_i}\otimes\Vect{\Sigma_j}\ar[r]^-{\mathcal{P}^k\otimes\mathcal{P}^k}\ar[d] & \kwr{i}{k}\otimes\kwr{j}{k}\ar[d] \\
\Vect_{\Sigma_i\times\Sigma_j}(X^{i+j})\ar[d] \ar[r]^{\mathcal{P}^k} & K_{_{(\Sigma_i\times\Sigma_j)\wr\Sigma_k}}(X^{i+j}) \ar[d]\\
\Vect_{\Sigma_{i+j}}\ar[r]^{\mathcal{P}^k} &\kwr{i+j}{k}}
$$
The top square is the fact that $(E\boxtimes F)^{\boxtimes k}\rightarrow E^{\boxtimes k}\boxtimes F^{\boxtimes k}$ is a map of $(\Sigma_i\times\Sigma_j)\wr\Sigma_k$ bundles.  The bottom square is Lemma \ref{C1}.

\end{proof}

\item[5.] The following commutes
$$
\xymatrix{ \K(X)\ar[r]^-{\tau}\ar[d]_{\tau} &\K P\ar[d]^{\dot{\tau}} \\
\K P\ar[r]^-{\Box} &\K Q}
$$
A definition for $\dot{\tau}$ is given in Section \ref{LambdaandTau}.  We repeat the details here.  We treat $R(\Sigma_l)$ as the $K$-theory ring $K_{_{\Sigma_l}}(\pt)$.  Then we get a map
$$
\tau^k\otimes\mathcal{P}^k: \K(X)\times K_{_{\Sigma_l}}(\pt)\longrightarrow \K(X)\otimes R(\Sigma_k)\otimes K_{_{\Sigma_l\wr\Sigma_k}}(\pt).
$$
We follow this map by the map 
$$
\operatorname{Id}\otimes\Ind{\Sigma_k}{\Sigma_l\wr\Sigma_k}\otimes \operatorname{Id}:\K(X)\otimes R(\Sigma_k)\otimes K_{_{\Sigma_l\wr\Sigma_k}}(\pt)\rightarrow \K(X)\otimes K_{_{\Sigma_l\wr\Sigma_k}}\otimes K_{_{\Sigma_l\wr\Sigma_k}}
$$ 
Then take the map $\operatorname{Id}_{\K(X)}\otimes\cdot$ where $\cdot$ is the internal tensor product in the ring $K_{_{\Sigma_l\wr\Sigma_k}}(\pt)$.  We define the map $\dot{\tau}^k_l$ to be the composite of all these maps 
$$
\dot{\tau}^k_l: \K(X)\times K_{_{\Sigma_l}}(\pt)\longrightarrow \K(X)\otimes K_{_{\Sigma_l\wr\Sigma_k}}(\pt).
$$
We define $\dot{\tau}_l$ to be the product over all $k$.  Then $\dot{\tau}_l$ is bi-additive from the additions to the cross product.  These maps extend to $\K(X)\otimes K_{_{\Sigma_l}}(\pt)$ and we define $\dot{\tau}$ to be the product of these maps over all $l$.  

The map $\Box$ is induced by the maps $R(\Sigma_{kl})\rightarrow R(\Sigma_l\wr\Sigma_k)$.

\begin{proof}
We divide the diagram up and show that the pieces commute:

\begin{center}
$$
\xymatrixrowsep{1cm}
\xymatrixcolsep{.1cm}
\xymatrix{ \ar @{} [dr] |{(1)} \ksig{m} \ar[r]^-{\mathcal{P}^l}\ar[d]_{\mathcal{P}^{kl}} & \ar @{} [dr] |{(2)} \kwr{m}{l} \ar[d]^{\mathcal{P}^k}\ar[r] & K_{_{\Sigma_{lm}\times\Sigma_l}}(X^{lm}) \ar[d]^{\mathcal{P}^k}\ar[r] & \ksig{lm}\otimes R(\Sigma_l) \ar[ddd]^{\dot{\tau}^k_l}\\
\kwr{m}{kl}\ar[d] \ar[r] & \ar @{} [d] |{(3)} K_{_{(\Sigma_m\wr\Sigma_l)\wr\Sigma_k}}(X^{klm})\ar[r] & \ar @{} [ddr] |{(5)} K_{_{(\Sigma_{lm}\times\Sigma_l)\wr\Sigma_k}}(X^{klm}) \ar[d]\\
K_{_{\Sigma_{klm}\times\Sigma_{kl}}}(X^{klm})\ar[d] \ar[rr] & & \ar @{} [dl] |{(4)} K_{_{\Sigma_{klm}\times (\Sigma_l\wr\Sigma_k)}}(X^{klm})\ar[dr] \\ 
\ksig{klm}\otimes R(\Sigma_{kl}) \ar[rrr] &  & & \ksig{klm}\otimes R(\Sigma_l\wr\Sigma_k) }     
$$
\end{center}

Square (1) commutes by Lemma \ref{L2}; (2) commutes Lemma \ref{C1}; and (4) is trivial.  Square (3) is another application of Mackey's Lemma using
$$
(\Sigma_m\wr\Sigma_{kl})\cap\Sigma_{klm}\times (\Sigma_l\wr\Sigma_k)=(\Sigma_m\wr\Sigma_l)\wr\Sigma_k
$$
embedded in $\Sigma_{klm}\times\Sigma_{kl}$.  

Square (5) commutes as follows.  Note first that both maps are homomorphic from + to the double cross product, when we product over all $k$.  For a proof of the right side the reader is referred to \cite[page 137]{Hoffman}.  The left side follows from the fact that $\prod_k\mathcal{P}^k(x+y)=\prod_k\mathcal{P}^k(x)\times\prod_k\mathcal{P}^k(y)$ and the diagram:
$$
\xymatrix{ K_{_{(\Sigma_{lm}\times\Sigma_l)\wr\Sigma_i}}(X^{lmi})\times K_{_{(\Sigma_{lm}\times\Sigma_l)\wr\Sigma_j}}(X^{lmj}) \ar[r] \ar[d] & K_{_{(\Sigma_{lm}\times\Sigma_l)\wr\Sigma_{i+j}}}(X^{(i+j)lm}) \ar[d] \\
K_{_{(\Sigma_l\wr\Sigma_i)\times\Sigma_{lmi}}}(X^{lmi})\times K_{_{(\Sigma_k\wr\Sigma_j)\times\Sigma_{lmj}}}(X^{lmj})\ar[r] & K_{_{\Sigma_{(i+j)lm}\times (\Sigma_l\wr\Sigma_{i+j})}}(X^{(i+j)lm}) }
$$
Therefore it suffices to check commutativity for $\ksig{lm}\times R(\Sigma_l)\subseteq K_{_{\Sigma_{lm}\times\Sigma_l}}(X^{lm})$.  The appropriate diagram is now:
$$
\xymatrixrowsep{1cm}
\xymatrixcolsep{.4cm}
\xymatrix{ \ksig{lm}\times R(\Sigma_l)\ar[r]_-{\mathcal{P}^k\times\mathcal{P}^k} \ar@{^{(}->}[dd] & \kwr{lm}{k}\times\kwr{l}{k} \ar@{^{(}->}[d] \ar[r] & K_{_{\Sigma_{klm}\times\Sigma_k}}(X^{klm})\times\kwr{l}{k} \ar@{^{(}->}[d]\\
& K_{_{(\Sigma_{lm}\wr\Sigma_k)\times (\Sigma_l\wr\Sigma_k)}}(X^{klm}) \ar[d] & K_{_{\Sigma_{klm}\times\Sigma_k\times (\Sigma_l\wr\Sigma_k)}}(X^{klm})\ar[d]\\
K_{_{\Sigma_{lm}\times\Sigma_l}}(X^{lm})\ar[r]^-{\mathcal{P}^k} & K_{_{(\Sigma_{lm}\times\Sigma_l)\wr\Sigma_k}}(X^{klm}) \ar[r] & K_{_{\Sigma_{klm}\times (\Sigma_l\wr\Sigma_k)}}(X^{klm}) }
$$
The right hand commutes by another application of Mackey's Lemma, this time with
$$
((\Sigma_{lm}\wr\Sigma_k)\times (\Sigma_l\wr\Sigma_k))\cap (\Sigma_{klm}\times (\Sigma_l\wr\Sigma_k))=(\Sigma_{lm}\times\Sigma_l)\wr\Sigma_k.
$$
The left hand commutes due to the fact that $(E\boxtimes F)^{\boxtimes k}\cong E^{\boxtimes k}\boxtimes F^{\boxtimes k}$ is an isomorphism of $(\Sigma_{lm}\times\Sigma_l)\wr\Sigma_k$-bundles.

\end{proof}

\end{itemize}

\section{A Thom Isomorphism}\label{ThomIso}

Now that most of our structure is in place, we can provide some results.  The first of these will be a Thom Isomorphism Theorem. 

In what follows it will be more convenient to use a different definition of $K_G(X,A)$, in keeping with definitions given in \cite{Atiyah}, \cite{Segal}.  Here, we let $X$ be any compact, Hausdorf space, with a closed subspace $A$.  By a complex of $G$-vector bundles over the pair $(X,A)$, we shall mean a complex of vector bundles, $E_i$, over $X$:
$$
\cdots 0\longrightarrow E_i\stackrel{d}{\longrightarrow} E_{i+1}\stackrel{d}{\longrightarrow} \cdots \stackrel{d}{\longrightarrow} E_{i+k}\longrightarrow 0\longrightarrow \cdots
$$
such that $d^2=0$.  We also require the complex to be acyclic when restricted to $A$.  A morphism of two complexes, $E_*$ and $F_*$, is a sequence of functions, $f_i:E_i\rightarrow F_i$, that commute with the differential $d$.  As with the other definition of the $K$ groups, the isomorphism classes of complexes form an abelian semigroup.  Two isomorphism classes, $E_*$ and $F_*$, are homotopic if there is complex of $G$-vector bundles over $(X\times [0,1],A\times [0,1])$ such that the restriction to $X\times\{0\}$ is $E_*$ and the restriction to $X\times\{1\}$ is $F_*$.  Furthermore, two homotopy classes of complexes, $E_*$ and $F_*$, are said to be equivalent if there are two complexes, $E_*'$ and $F_*'$, acyclic over $X$, such that
$$
E_*\oplus E_*' \simeq F_*\oplus F_*'.
$$
The set of equivalence classes can be shown to be isomorphic to $K_G(X,A)$, \cite[page 148]{Segal}.  The isomorphism is simply $E_*\mapsto\sum_k(-1)^kE_k$. 

In what follows we will assume that $A=\emptyset$.  At the end of the section, we will discuss how to prove the result for arbitrary $A$. 

Given a $\Sigma_n$ complex, $E_*$, over $X^n$ and a $\Sigma_m$ complex, $F_*$, over $X^m$, define the product of $E_*$ and $F_*$ over $X^{n+m}$ to have $k$th term
$$
(E_*F_*)_k=\bigoplus_{p+q=k}\Ind{\Sigma_n\times\Sigma_m}{\Sigma_{n+m}}(E_p\boxtimes F_q).
$$
The $\tau$-ring structure, $\lambda$-ring structure, and comultiplication (when appropriate) are defined similarly.

Recall that if $E\rightarrow X$ is a $G$-vector bundle then we have an additive isomorphism 
$$
T_E:K_G(X)\longrightarrow K_G(D(E),S(E))
$$ 
called the Thom isomorphism.  Here $D(E)$ and $S(E)$ are the associated disk and sphere bundles.  We need to recall precisely how this is defined.

Let the bundle map be $p:E\rightarrow X$.  Then the pullback bundle $p^*E\rightarrow E$ has a natural section along the diagonal $\delta:E\rightarrow E\times_XE$.  Now we form the Koszul complex 
$$
\Lambda^E_*=\cdots \longrightarrow 0\longrightarrow \C\stackrel{d}{\longrightarrow}\Lambda^1p^*E\stackrel{d}{\longrightarrow}\Lambda^2p^*E\stackrel{d}{\longrightarrow}\cdots,
$$
where $d(v)=v\wedge\delta(x)$ for $v\in\Lambda^ip^*E_x$.  For our purposes, we must restrict each $\Lambda^ip^*E$ to $D(E)$, but we will keep the same notation.  Since $\delta$ is non-vanishing on $S(E)$ this complex will be acyclic when restricted to the sphere bundle, \cite[page 99]{Atiyah}.  If $F_*$ is a complex on $X$ then $p^*F_*$ is a complex on $E$ (again restrict down to $D(E)$) and the Thom homomorphism is then
$$
T_E(F_*)=\Lambda^E_*\otimes p^*F_*.
$$
Note that, if $E$ is a bundle over $X$, then $E^n$ is a $\Sigma_n$-bundle over $X^n$.  Therefore $T_{E^n}$ is an isomorphism
$$
K_{_{\Sigma_n}}(X^n)\longrightarrow K_{_{\Sigma_n}}(D(E^n),S(E^n)).
$$
But, for each $n$, $D(E^n)\cong D(E)^n$ and  
$$
[D(E)^{n-1}\times S(E)]\cup \cdots \cup [S(E)\times D(E)^{n-1}]\cong S(E^n).
$$
This gives a homeomorphism of pairs $(D(E^n),S(E^n))\cong (D(E),S(E))^n$.  With a slight abuse of notation, we have isomorphisms for all $n$:
$$
T_{E^n}:K_{_{\Sigma_n}}(X^n)\stackrel{\cong}{\longrightarrow} K_{_{\Sigma_n}}((D(E),S(E))^n).
$$
Adding all these together, we get an isomorphism of abelian groups
$$
T:\K(X)\longrightarrow \K(D(E),S(E)).
$$

We now wish to show that this isomorphism agrees with both the ring and the $\lambda$-ring structure.

For the ring structure we begin with two elements $V_*\in\ksig{n}$ and $W_*\in\ksig{m}$.  If we do the Thom homomorphism to each separately and then multiply we get:
\begin{align*}
\mathrlap{\Ind{\Sigma_n\times\Sigma_m}{\Sigma_{n+m}}\left[\bigoplus_{p+q=k}\left(\bigoplus_{i+j=p,\, s+t=q}\left((\Lambda^ip^*E^n\otimes V_j)\boxtimes (\Lambda^sp^*E^m\otimes W_t)\right)\right)\right]} \\
\hphantom{\qquad} & \\
&\cong \bigoplus_{i+j+s+t=k}\Ind{\Sigma_n\times\Sigma_m}{\Sigma_{n+m}}\left[ (\Lambda^ip^*E^n\otimes V_j)\boxtimes (\Lambda^sp^*E^m\otimes W_t)\right] \\
&\cong \bigoplus_{i+j+s+t=k}\Ind{\Sigma_n\times\Sigma_m}{\Sigma_{n+m}}\left[ (\Lambda^ip^*E^n\boxtimes\Lambda^sp^*E^m)\otimes (V_j\boxtimes W_t)\right] \\
&\cong\bigoplus_{t+j+p=k}\Ind{\Sigma_n\times\Sigma_m}{\Sigma_{n+m}}\left[\left(\bigoplus_{i+s=p}\Lambda^ip^*E^n\boxtimes\Lambda^sp^*E^m\right)\otimes (V_j\boxtimes W_t)\right] \\
&\cong \bigoplus_{j+t+p=k}\Ind{\Sigma_n\times\Sigma_m}{\Sigma_{n+m}}\left[\Res{\Sigma_n\times\Sigma_m}{\Sigma_{n+m}}\left(\Lambda^pp^*E^{n+m}\right)\otimes (V_j\boxtimes W_t)\right] \\
&\cong \bigoplus_{j+t+p=k}\Lambda^pp^*E^{n+m}\otimes\Ind{\Sigma_n\times\Sigma_m}{\Sigma_{n+m}}(V_j\boxtimes W_t),
\end{align*}
as $T(V_*W_*)_k$.  The second isomorphism is due to the fact that
$$
(E_1\otimes E_2)\boxtimes (F_1\otimes F_2)\cong (E_1\boxtimes F_1)\otimes (E_2\boxtimes F_2)
$$
as $\Sigma_n\times\Sigma_m$ bundles.  The fourth isomorphism is due jointly to the facts that 
$$
p^*E^{n+m}\cong (p^*E)^{n+m}
$$ 
and 
$$
\Lambda^pF^{n+m}\cong \bigoplus_{i+j=p}\Lambda^iF^n\boxtimes \Lambda^jF^m
$$
as $\Sigma_n\times\Sigma_m$ bundles.  The final isomorphism is an instance of Frobenius reciprocity:
$$
\Ind{H}{G}\left[\Res{H}{G}(E_1)\otimes E_2\right]\cong E_1\otimes\Ind{H}{G}(E_2)
$$
for a $G$ bundle $E_1$ and an $H$ bundle $E_2$, $H$ a subgroup of a group $G$.

Since 
$$
\bigoplus_{j+t+p=k}\Lambda^pp^*E^{n+m}\otimes\Ind{\Sigma_n\times\Sigma_m}{\Sigma_{n+m}}(V_j\boxtimes W_t)
$$
is what we get by multiplying and then applying the Thom homomorphism, we are done.  

We now move onto the $\lambda$-ring isomorphism.  If we start by taking $V_*\in\ksig{n}$, apply $\lambda^m$, and then do the Thom homomorphism we get:
\begin{align*}
\mathrlap{\bigoplus_{p+q=k}\Lambda^pp^*E^{nm}\otimes\Ind{\Sigma_n\wr\Sigma_m}{\Sigma_{nm}}\left(\sgn_m\otimes\bigoplus_{i_1+i_2+\cdots +i_m=q}V_{i_1}\boxtimes V_{i_2}\boxtimes\cdots\boxtimes V_{i_m}\right)} \\
\hphantom{\qquad} & \\
&\cong \Ind{\Sigma_n\wr\Sigma_m}{\Sigma_{nm}}\left(\sgn_m\otimes\bigoplus_{i_1+\cdots +i_m+j_1+\cdots +j_m=k}\Lambda^{i_1}p^*E^n\otimes V_{j_1}\boxtimes\cdots\boxtimes\Lambda^{i_m}p^*E^n\otimes V_{j_m}\right)
\end{align*}
as $(T(\lambda^m(V_*)))_k$.  Since the second line is precisely what we get when we do $T$ then $\lambda^m$, we are done.  Thus $T$ is an isomorphism as $\lambda$-rings. 

We observe that the discussion generalizes to the pair $(X,A)$.  So that our work will fit on the page, let us denote
$$
A(n)=(X^{n-1}\times A)\cup (X^{n-2}\times A \times X)\cup\ldots\cup (A\times X^{n-1}).
$$
At level $n$ we have a Thom isomorphism
$$
T_{E^n}:K_{_{\Sigma_n}}(X^n,A(n))\longrightarrow K_{_{\Sigma_n}}(D(E^n),S(E^n)\cup D(E^n)|_{A(n)}). 
$$
When $n=1$, this is
$$
K(D(E),S(E)\cup D(E)|_A).
$$
As before we need that $D(E^n)\cong D(E)^n$ and $S(E^n)\cong S(E)^n$.  What we also need is that
$$
[D(E)^{n-1}\times (S(E)\cup D(E)|_A)]\cup\ldots \cup [(S(E)\cup D(E)|_A)\times D(E)^{n-1}]\cong S(E)^n\cup D(E)^n|_{A(n)}.
$$
But, this is also true:  the collection
$$
(D(E)^{n-1}\times S(E))\cup\ldots\cup (S(E)\times D(E)^{n-1}
$$
will give us $S(E)^n$.  Now all we would need is for 
$$
(D(E)^{n-1}\times D(E)|_A)\cup\ldots\cup (D(E)|_A\times D(E)^{n-1})\cong D(E)^n|_{A(n)}.
$$
The left hand side is readily seen to be the same as restricting $D(E)^n$ down to
$$
A(n)=(X^{n-1}\times A)\cup\ldots\cup (A\times X^{n-1}).
$$
The proof of the Thom isomorphism above with $A=\emptyset$ generalizes to give
$$
\K(X,A)\stackrel{\cong}{\longrightarrow} \K(D(E),S(E)\cup D(E)|_A).
$$

\section{Spaces with Even Cells}\label{FreeIso}

By induction on cells it is not hard to see that if $X$ is a CW-complex made only from finitely many, even dimensional cells, then $K(X)$ is the free abelian group on generators in correspondence with those cells and $K^1(X)=0$.  What will be shown here is that, if $X$ has only even dimensional cells, then $\K(X)$ is the free $\lambda$-ring on $K(X)$ as defined in Section \ref{Free}.

If $X$ is a point, then we already know that $\K(X)\cong\oplus R\Sigma_n$, which is the free $\lambda$-ring on one generator.  

Let $X$ be a space constructed out of even cells.  Assume $\K(X)$ is isomorphic to $\mathfrak{F}(K(X))$.  We assume that $Y$ is obtained from $X$ by adding one even dimensional cell, denoted by $D$, via attaching along a sphere $S$, and proceed by induction.  The strategy of the proof will be to first show that $\K(Y)$ is a polynomial ring.  We can then more easily show that the induced map $\mathfrak{F}(K(Y))\rightarrow \K(Y)$ is an isomorphism.  We shall also prove along the way that $\K^1(Y)=0$, where we define
$$
\K^1(Y,X)=\bigoplus_{n\geq 0}K^1_{_{\Sigma_n}}((Y,X)^n).
$$
We then have the additional assumption that $\K^1(X)=0$.  Again we know this to be true for $X=\operatorname{pt}$:  In this case 
$$
K^1_{_{\Sigma_n}}(\operatorname{pt})\cong K^1(\operatorname{pt})\otimes R(\Sigma_n)=0
$$
To avoid confusion we shall begin writing $\K^0(X)$ for $\K(X)$.

Now, $Y^n$ has a filtration:
$$
\emptyset = F^n_{-1}\subseteq F_0^n=X^n\subseteq\cdots\subseteq F^n_{n-1}\subseteq F^n_n=Y^n
$$
where $F^n_k$ is the set of points in $Y^n$ such that at least $n-k$ of the coordinates are in $X$, or at most $k$ of them are in $Y$ but not $X$.   This leads to a nested set of pairs:
$$
(Y^n,F^n_{-1})\subseteq (Y^n,F^n_0)\subseteq \cdots\subseteq (Y^n,F^n_{n-1})\subseteq (Y^n,F^n_n)
$$
which gives rise to a filtration of $K^*_{_{\Sigma_n}}(Y^n)$:
$$
K^*_{_{\Sigma_n}}(Y^n,F^n_{-1})    \supseteq  K^*_{_{\Sigma_n}}(Y^n,F^n_0)\supseteq\cdots\supseteq K^*_{_{\Sigma_n}}(Y^n,F^n_{n-1})\supseteq K^*_{_{\Sigma_n}}(Y^n,Y^n).
$$
It might also be good to note that this filtration is 
$$
K^*_{_{\Sigma_n}}(Y^n) \supseteq K^*_{_{\Sigma_n}}(Y^n,X^n)\supseteq \cdots\supseteq K^*_{_{\Sigma_n}}((Y,X)^n) \supseteq K^*_{_{\Sigma_n}}(Y^n,Y^n).
$$
The resulting abelian groups of the spectral sequence are:
$$
E^{n,m,*}_1=\frac{K^*_{_{\Sigma_{n+m}}}(Y^{n+m},F^{n+m}_n)}{K^*_{_{\Sigma_{n+m}}}(Y^{n+m},F^{n+m}_{n+1})}.
$$

The following lemma shows that this is actually a sequence of commutative algebras.

\begin{lemma}\label{L3}
The induction product gives a map
$$
K^*_{_{\Sigma_n}}(Y^n,F^n_k)\otimes K^*_{_{\Sigma_m}}(Y^m,F^m_l)    \longrightarrow K^*_{_{\Sigma_{n+m}}}(Y^{n+m},F^{n+m}_{k+l}).
$$
\begin{proof}
Take $K^*_{_{\Sigma_n}}(Y^n,F^n_k)$ and $K^*_{_{\Sigma_m}}(Y^m,F^m_l)$.  Then the external tensor product goes
$$
K^*_{_{\Sigma_n}}(Y^n,F^n_k)\otimes K^*_{_{\Sigma_m}}(Y^m,F^m_l)\longrightarrow K^*_{_{\Sigma_n\times\Sigma_m}}(Y^{n+m},(Y^n\times F^m_l)\cup (F^n_k\times Y^m)).
$$
But, 
$$
(Y^n\times F^m_l)\cup (F^n_k\times Y^m)\supseteq F^{n+m}_{k+l}.
$$ 
Thus we get a map
$$
K^*_{_{\Sigma_n\times\Sigma_m}}(Y^{n+m},(Y^n\times F^m_l)\cup (F^n_k\times Y^m))\longrightarrow K^*_{_{\Sigma_n\times\Sigma_m}}(Y^{n+m},F^{n+m}_{k+l}).
$$
The induction map then lands in 
$$
K^*_{_{\Sigma_{n+m}}}(Y^{n+m},F^{n+m}_{k+l}).
$$  
\end{proof}
\end{lemma}

So, we now have a spectral sequence of commutative algebras.  Since 
$$
K^1_{_{\Sigma_m}}(Y^m,F^m_n)=K^1_{_{\Sigma_m}}(F^m_n)=0,
$$
the long exact sequence of the pair gives a map of short exact sequences
$$
\xymatrix{       
 0 \ar[r] & K^*_{_{\Sigma_m}}(Y^m,F^m_{n+1}) \ar[r]\ar[d]  & K^*_{_{\Sigma_m}}(Y^m)  \ar[r]\ar[d]^{\operatorname{Id}} & K^*_{_{\Sigma_m}}(F^m_{n+1})  \ar[r]\ar[d] & 0 \\
0 \ar[r] & K^*_{_{\Sigma_m}}(Y^m,F^m_n) \ar[r] &  K^*_{_{\Sigma_m}}(Y^m) \ar[r]  &  K^*_{_{\Sigma_m}}(F^m_n)  \ar[r]  & 0     }
$$
The Snake Lemma then gives us a sequence:
$$
0\longrightarrow K^*_{_{\Sigma_m}}(F^m_{n+1},F^m_n)\stackrel{\cong}{\longrightarrow} E^{n,m-n,*}_1 \longrightarrow 0.
$$
Therefore
$$
E^{n,m,*}_1\cong K^*_{_{\Sigma_{n+m}}}(F^{n+m}_{n+1},F^{n+m}_n).
$$
Our goal is to show that $E^{*,*,0}$ is polynomial and $E^{*,*,1}=0$.  Then we will know that $\K^0(Y)$ is polynomial and $\K^1(Y)=0$.  We do this by proving $\K^*(Y,X)\otimes\K^*(X)\cong E_1$ and then show that $\K^*(Y,X)$ is polynomial with $\K^1(Y,X)=0$.  The ring $\K^0(X)$ is polynomial with $\K^1(X)=0$ by hypothesis.  First we show:

\begin{lemma}\label{L4}
The external tensor product
$$
K^*_{_{\Sigma_n}}((Y,X)^n)\otimes K^*_{_{\Sigma_m}}(X^m)\longrightarrow K^*_{_{\Sigma_{n+m}}}(F^{n+m}_n,F^{n+m}_{n-1})
$$
is an isomorphism.

\begin{proof}
We have maps of $\Sigma_{n+m}$-spaces 
$$
\Sigma_{n+m}\times_{_{\Sigma_n\times\Sigma_m}}(F^n_{n-1}\times X^m)\longrightarrow F^{n+m}_{n-1}
$$ 
and 
$$
\Sigma_{n+m}\times_{_{\Sigma_n\times\Sigma_m}}(Y^n\times X^m)\longrightarrow F^{n+m}_n
$$
that fit into a pushout square
$$
\xymatrix{ \Sigma_{n+m}\times_{_{\Sigma_n\times\Sigma_m}}(F^n_{n-1}\times X^m)\ar[r]\ar[d] & F^{n+m}_{n-1}\ar[d]\\
\Sigma_{n+m}\times_{_{\Sigma_n\times\Sigma_m}}(Y^n\times X^m) \ar[r] & F^{n+m}_n}
$$
Therefore, 
\begin{align*}
K^*_{_{\Sigma_{n+m}}}(F^{n+m}_n,F^{n+m}_{n-1}) & \cong  K^*_{_{\Sigma_{n+m}}}(\Sigma_{n+m}\times_{_{\Sigma_n\times\Sigma_m}}(Y^n\times X^m),\Sigma_{n+m}\times_{_{\Sigma_n\times\Sigma_m}}(F^n_{n-1}\times X^m))\\
&\cong  K^*_{_{\Sigma_n\times\Sigma_m}}((Y^n\times X^m),(F^n_{n-1}\times X^m))\\
 &\cong   K^*_{_{\Sigma_n}}((Y,X)^n)\otimes K^*_{_{\Sigma_m}}(X^m)
\end{align*}
The third isomorphism is given by the external product:
$$
K_{_{\Sigma_n}}((Y,X)^n)\otimes K_{_{\Sigma_m}}(X^m)\longrightarrow K_{_{\Sigma_n\times\Sigma_m}}(Y^n\times X^m,F^n_{n-1}\times X^m).
$$
The Thom isomorphism tells us that $K_{_{\Sigma_n}}((Y,X)^n)$ is finitely generated and free as an abelian group.  So, the external product is an isomorphism by the Kunneth formula.

\end{proof}
\end{lemma}

We have now shown that $E^{n,m,1}=0$, hence $\K^1(Y)=0$.  Due to this we will stop using $\K^0$ and go back to just $\K$.  Since $K_{_{\Sigma_{n+m}}}(F^{n+m}_n,F^{n+m}_{n-1})\cong E^{n-1,m+1}_1$, we have shown that the map 
$$
\K(Y,X)\otimes \K(X)\longrightarrow E_1
$$ 
is an isomorphism of bigraded abelian groups.  It remains to show that it is an isomorphism of bigraded commutative algebras.  The next lemma demonstrates that.
\begin{lemma}\label{L5}
The following diagram commutes:
$$
\xymatrix{    K_{_{\Sigma_n}}((Y,X)^n)\otimes K_{_{\Sigma_m}}(X^m) \ar[r]^-\cong & K_{_{\Sigma_{n+m}}}(F^{n+m}_n,F^{n+m}_{n-1}) \\
K_{_{\Sigma_n}}((Y,X)^n)\otimes K_{_{\Sigma_m}}(Y^m) \ar[u] \ar[r] \ar[d]&   K_{_{\Sigma_{n+m}}}(Y^{n+m},F^{n+m}_{n-1}) \ar[u] \ar[d]  \\
K_{_{\Sigma_n}}(Y^n)\otimes K_{_{\Sigma_m}}(Y^m) \ar[r] &  K_{_{\Sigma_{n+m}}}(Y^{n+m})       }
$$
where the vertical maps are induced by inclusion, the top map is the isomorphism from Lemma \ref{L4}, the bottom map is the ring multiplication, and the middle map is external tensor
$$
K_{_{\Sigma_n}}((Y,X)^n)\otimes K_{_{\Sigma_m}}(Y^m)      \stackrel{\boxtimes}{\longrightarrow}      K_{_{\Sigma_n\times\Sigma_m}}(Y^{n+m},(F^n_{n-1}\times Y^m))   
$$
followed by
$$
 K_{_{\Sigma_n\times\Sigma_m}}(Y^{n+m},(F^n_{n-1}\times Y^m))   \subseteq    K_{_{\Sigma_n\times\Sigma_m}}(Y^{n+m},F^{n+m}_{n-1})   \stackrel{\operatorname{Ind}}{\longrightarrow}          K_{_{\Sigma_{n+m}}}(Y^{n+m},F^{n+m}_{n-1}) .
$$
\begin{proof}
We need the fact that $(Y,X)^n=(Y^n,F^n_{n-1})$.  To show that the top square commutes, break it up into the following:
$$
\xymatrix{      K_{_{\Sigma_n}}(Y^n,F^n_{n-1})\otimes K_{_{\Sigma_m}}(Y^m) \ar[r] \ar[d]^-\boxtimes &   K_{_{\Sigma_n}}(Y^n,F^n_{n-1})\otimes K_{_{\Sigma_m}}(X^m)\ar[d]^-\boxtimes   \\
 K_{_{\Sigma_n\times\Sigma_m}}(Y^{n+m},F^n_{n-1}\times Y^m) \ar[d] \ar[r]   &   K_{_{\Sigma_n\times\Sigma_m}}(Y^n\times X^m,F^n_{n-1}\times X^m) \ar[d]^\cong    \\
 K_{_{\Sigma_n\times\Sigma_m}}(Y^{n+m},F^{n+m}_{n-1})  \ar[d]^{\operatorname{Ind}} \ar[r]   &      K_{_{\Sigma_{n+m}}}(\Sigma_{n+m}\times_{_{\Sigma_n\times\Sigma_m}}(Y^n\times X^m,F^n_{n-1}\times X^m)) \ar[d]^\cong    \\
 K_{_{\Sigma_{n+m}}}(Y^{n+m},F^{n+m}_{n-1}) \ar[r]^-f   &      K_{_{\Sigma_{n+m}}}(F^{n+m}_n,F^{n+m}_{n-1})          }
$$
The top two squares commute trivially.  It takes more care to show that the bottom square commutes.  To show this we will show that the diagram
$$
\xymatrix{ 
K_{_{\Sigma_n\times\Sigma_m}}(Y^{n+m},F^{n+m}_{n-1})  \ar[d]^{\operatorname{Ind}} \ar[r]     &     K_{_{\Sigma_{n+m}}}(\Sigma_{n+m}\times_{_{\Sigma_n\times\Sigma_m}}(Y^n\times X^m,F^n_{n-1}\times X^m))     \\
 K_{_{\Sigma_{n+m}}}(Y^{n+m},F^{n+m}_{n-1}) \ar[r]^-f   &      K_{_{\Sigma_{n+m}}}(F^{n+m}_n,F^{n+m}_{n-1}) \ar[u]^\cong         }
$$
commutes where the map 
$$
K_{_{\Sigma_{n+m}}}(F^{n+m}_n,F^{n+m}_{n-1})  \longrightarrow   K_{_{\Sigma_{n+m}}}(\Sigma_{n+m}\times_{_{\Sigma_n\times\Sigma_m}}(Y^n\times X^m,F^n_{n-1}\times X^m))
$$
is the isomorphism from the pushout diagram in Lemma \ref{L4} which is inverse to the isomorphism in the diagram we now consider.  We start with a bundle in $K_{_{\Sigma_n\times\Sigma_m}}(Y^{n+m},F^{n+m}_{n-1})$ which is of the form $E\boxtimes F$ where $E$ restricted to $F^n_{n-1}$ is trivial.  We take this bundle and induce up $\Ind{\Sigma_n\times\Sigma_m}{\Sigma_{n+m}}(E\boxtimes F)$.  We then restrict this bundle down to $Y^n\times X^m$ and only consider the action of $\Sigma_n\times\Sigma_m$.  The induced bundle is a direct sum of copies of $E\boxtimes F$, 
one for each coset of $\Sigma_{n+m}/(\Sigma_n\times\Sigma_m)$.  Write this as
$$
\xi=(E\boxtimes F)_{\sigma_1}\oplus\cdots\oplus (E\boxtimes F)_{\sigma_r},
$$
where $\sigma_i$ are coset representatives.  If we have a point
$$
\bar{a}=(a_1,\dots,a_n,a_{n+1},\dots,a_{n+m})\in Y^n\times X^m
$$
we may write the preimage of $\bar{a}$ under the projection
$$
\Sigma_{n+m}\times_{\Sigma_n\times\Sigma_m}Y^{n+m}\longrightarrow Y^{n+m}
$$
as the collection of points
$$
\left(\sigma_i,a_{\sigma_i^{-1}(1)},\dots,a_{\sigma^{-1}_i(n+m)}\right)
$$
for $i=1,\dots,r$.  If $\sigma_i\notin \Sigma_n\times\Sigma_m$ then for some $k$, with $1\leq k\leq n$, we have $\sigma_i(k)>n$.  But this forces $E$ to be restricted to a subspace of $F^n_{n-1}$ over which it is trivial.  Therefore $\xi$ is isomorphic to $E\boxtimes F$ restricted to the subspace $Y^n\times X^m$ which is what we get when we go along the top of the square.

To complete the lemma we need only show that the following commutes:
$$
\xymatrix{      K_{_{\Sigma_n}}(Y^n,F^n_{n-1})\otimes K_{_{\Sigma_m}}(Y^m) \ar[r] \ar[d]^-\boxtimes &   K_{_{\Sigma_n}}(Y^n)\otimes K_{_{\Sigma_m}}(Y^m)\ar[d]^-\boxtimes   \\
 K_{_{\Sigma_n\times\Sigma_m}}(Y^{n+m},F^n_{n-1}\times Y^m) \ar[d] \ar[r]    &     K_{_{\Sigma_n\times\Sigma_m}}(Y^{n+m})  \ar[dd]^{\operatorname{Ind}}    \\
 K_{_{\Sigma_n\times\Sigma_m}}(Y^{n+m},F^{n+m}_{n-1})  \ar[d]^{\operatorname{Ind}}     &                 \\
 K_{_{\Sigma_{n+m}}}(Y^{n+m},F^{n+m}_{n-1}) \ar[r]   &      K_{_{\Sigma_{n+m}}}(Y^{n+m})     }
$$
But, this is just the fact that induction is natural.

\end{proof}

\end{lemma}

We are nearing the end.  Recall that $\mathfrak{F}(K(Y))$ is the free $\lambda$-ring on $K(Y)$.  We have the following commutative diagram:
$$
\xymatrix{     \mathfrak{F}(K(Y,X)) \ar[r]\ar[d]^-{\cong} & \mathfrak{F}(K(Y)) \ar[r]\ar[d] & \mathfrak{F}(K(X)) \ar[d]^-{\cong} \\
             \K(Y,X) \ar[r]  &  \K(Y) \ar[r]   &  \K(X)    }
$$
where the left map is an isomorphism induced by the map of pairs $(D,S)\longrightarrow (Y,X)$ and the right map is an isomorphism by assumption.  By definition $\mathfrak{F}(K(Y))$ is free as an algebra over $\mathfrak{F}(K(Y,X))$.  And $\K(Y)$ is free as an algebra over $\K(Y,X)$ by the work above.  The three lemmas of the section have demonstrated that we have an isomorphism
$$
\K(Y)\otimes_{\K(Y,X)}\mathbb{Z}\stackrel{\cong}{\longrightarrow} \K(X).
$$
This fits neatly into a diagram
$$
\xymatrix{        \mathfrak{F}(K(Y))\otimes_{\mathfrak{F}(K(Y,X))}\mathbb{Z} \ar[r]^-{\cong}\ar[d] & \mathfrak{F}(K(X)) \ar[d]^-{\cong}\\
                  \K(Y)\otimes_{\K(Y,X)} \mathbb{Z} \ar[r]^-{\cong} & \K(X)    }
$$
Therefore we have that
$$
 \mathfrak{F}(K(Y))\otimes_{\mathfrak{F}(K(Y,X))}\mathbb{Z}\longrightarrow  \K(Y)\otimes_{\K(Y,X)} \mathbb{Z}
$$
is an isomorphism.  Since $\mathfrak{F}(K(Y,X))\cong \K(Y,X)$ naturally, the isomorphism becomes
$$
\mathfrak{F}(K(Y))\otimes_{\mathfrak{F}(K(Y,X))}\mathbb{Z}\stackrel{\cong}{\longrightarrow} \K(Y)\otimes_{\mathfrak{F}(K(Y,X))} \mathbb{Z}.
$$
Because everything is a polynomial algebra, we may conclude that $\mathfrak{F}(K(Y))$ is isomorphic to $\K(Y)$.
\begin{theorem}\label{EvenFree}
If $X$ is a CW-complex, constructed from a finite number of even dimensional cells, then the natural map of $\lambda$-rings
$$
\mathfrak{F}(K(X))\longrightarrow \K(X)
$$
induces an isomorphism. 
\end{theorem}

Just as in the case of the Thom isomorphism, we can extend these results to pairs, $(X,A)$.  In this case, Theorem \ref{EvenFree} becomes:

\begin{theorem}\label{EvenFreePair}
If $X$ is a CW-complex, constructed by attaching a finite number of even dimensional cells to a finite CW-complex $A$, then the natural map of $\lambda$-rings
$$
\mathfrak{F}(K(X,A))\longrightarrow \K(X,A)
$$
induces an isomorphism.
\end{theorem}

The way to see this is that the base case $\K(A,A)$ is clearly true.  One then checks that the lemmas from this section can be appropriately modified.  So, if we know it to be true for the pair $(X,A)$, and $Y$ is built from $X$ by attaching an even dimensional cell, our new diagram becomes:
$$
\xymatrix{   \mathfrak{F}(K(Y,X)) \ar[r]\ar[d]^-{\cong} & \mathfrak{F}(K(Y,A)) \ar[r]\ar[d] & \mathfrak{F}(K(X,A)) \ar[d]^-{\cong} \\
             \K(Y,X) \ar[r]  &  \K(Y,A) \ar[r]   &  \K(X,A)    }
$$

\bibliography{Thesis}{}
\bibliographystyle{plain}
\end{document}